\documentclass{siamltex}
\usepackage{amsmath}
\usepackage{latexsym}

\def\noprint#1{}

\newtheorem{assumption}{Assumption}
\newtheorem{Exa}{Example}[section]

\newcommand{\gph}{\mbox{\rm gph}}

\newcommand{\ip}[2]{\mbox{$\langle #1,#2 \rangle$}}
\def\tto{\;{\lower 1pt \hbox{$\rightarrow$}}\kern -10pt
           \hbox{\raise 2.8pt \hbox{$\rightarrow$}}\;}

\newcommand{\beq}{\begin{equation}}
\newcommand{\eeq}{\end{equation}}
\newcommand{\IR}{\makebox{\sf I \hspace{-8.5 pt} R \hspace{-4.0pt}}}
\newcommand{\R}{\IR}

\newcommand{\cC}{{\cal C}}

\def\eqnok#1{(\ref{#1})}

\newcommand{\dist}{\mbox{\rm dist}}

\def\noprint#1{}

\newcommand{\btab}{\begin{tabbing}
\ \ \= thenn \= thenn \= thenn \= thenn \= thenn \= \kill}
\newcommand{\etab}{\end{tabbing}}

\title{Identifying Activity\footnote{\today}}

\author{A.~S. Lewis\thanks{ORIE, Cornell University, Ithaca, NY 14853, U.S.A.
\texttt{people.orie.cornell.edu/\~~\hspace{-4pt}aslewis aslewis\char64 orie.cornell.edu}.
Research supported in part by National Science Foundation Grant DMS-0806057.} 
\and 
S.~J. Wright\thanks{Computer Sciences Department, University of Wisconsin, 1210 W.~Dayton Street, Madison, WI 53706.
\texttt{swright\char64 cs.wisc.edu 
pages.cs.wisc.edu/\~~\hspace{-4pt}swright}.
Research supported in part by National Science Foundation Grant 0430504.}}

\begin{document}

\maketitle

\begin{abstract}
  Identification of active constraints in constrained optimization is
  of interest from both practical and theoretical viewpoints, as it
  holds the promise of reducing an inequality-constrained problem to
  an equality-constrained problem, in a neighborhood of a solution. We
  study this issue in the more general setting of composite nonsmooth
  minimization, in which the objective is a composition of a smooth
  vector function $c$ with a lower semicontinuous
  function $h$, typically nonsmooth but structured. In this setting, the graph of the generalized gradient $\partial h$ can often be decomposed into a union (nondisjoint) of simpler subsets.  ``Identification'' amounts to deciding which subsets of the graph are ``active'' in the criticality conditions at a given solution.  We give conditions under which any convergent sequence of approximate critical points finitely identifies the activity.
Prominent among these properties is a condition akin to
  the Mangasarian-Fromovitz constraint qualification, which ensures
  boundedness of the set of multiplier vectors that satisfy the optimality
  conditions at the solution. 
\end{abstract}



\section{Introduction} \label{sec:intro}


We study ``active set'' ideas for a composite optimization problem of the form
\begin{equation} \label{hc}
\min_x h \big( c(x) \big).
\end{equation}
Throughout this work, we make the following rather standard blanket assumption.
\medskip

\begin{assumption} \label{blanket} The function $h \colon \R^m \to
  \bar\R$ is lower semicontinuous and the function $c \colon \R^n \to
  \R^m$ continuously differentiable.  The point $\bar x \in \R^m$ is
  critical for the composite function $h \circ c$, and satisfies the
  condition
\begin{equation} \label{cq}
\partial^{\infty} h \big( c(\bar x) \big) \cap N(\nabla c(\bar x)^*) ~=~ \{0\}.
\end{equation}
\end{assumption}

\vspace{-12pt}
\noindent
In condition (\ref{cq}), $N(\cdot)$ denotes the null space and
$\partial^{\infty}$ denotes the horizon subdifferential, defined
below.

Some comments are in order.  Because the outer function $h$ can take
values in the extended reals $\bar\R = [-\infty,+\infty]$, we can
easily model constraints.  In many typical examples, $h$ is convex.
We develop the general case, although noting throughout how the theory
simplifies in the convex case.  For notational simplicity, we suppose
that the inner function $c$ is everywhere defined, the case where its
domain is an open subset being a trivial extension.  By a {\em
  critical} point for $h \circ c$, we mean a point satisfying the
condition $0 \in \partial (h \circ c)(\bar x)$.  Here, $\partial$
denotes the subdifferential of a nonsmooth function.  We refer to the
monographs \cite{Cla98,Roc98,Mor06} for standard ideas from
variational analysis and nonsmooth optimization, and in particular we
follow the notation and terminology of \cite{Roc98}.  For continuously
differentiable functions, the subdifferential coincides with the
derivative, while for convex functions it coincides with the classical
convex subdifferential.  Equation \eqnok{cq} is called a regularity
(or {\em transversality\/}) condition: $\partial^{\infty}$ denotes the
{\em horizon subdifferential}.  If the function $h$ is lower
semicontinuous, convex, and finite at the point $\bar c$, then
$\partial^{\infty} h(\bar c)$ is the normal cone (in the sense of
classical convex analysis) to the domain of $h$ at $\bar c$.  If in
addition $h$ is continuous at $\bar c$, then we have
$\partial^{\infty} h(\bar c) = \{0\}$.

A standard chain rule ensures the inclusion
\[
\partial ( h \circ c)(\bar x) ~\subset~ \nabla c(\bar x)^* \partial h \big( c(\bar x) \big).
\]
We deduce that there exists a vector $v \in \R^m$ satisfying the
conditions
\begin{equation} \label{eq:vcrit}
v \in \partial h \big( c(\bar x) \big), ~~ \nabla c(\bar x)^* v = 0.
\end{equation}
By analogy with classical nonlinear programming (as we shall see), we make the following definition.
\begin{definition}
A vector $v \in \R^m$ satisfying the conditions \eqnok{eq:vcrit} is called a {\bf multiplier vector} for the critical point $\bar x$.
\end{definition}

In seeking to solve the problem (\ref{hc}), we thus look
for a pair $(x,v) \in \R^n \times \R^m$ such that
\begin{equation} \label{xv.crit}
v \in \partial h \big( c(x) \big), ~~ \nabla c(x)^* v = 0.
\end{equation}
As we have just observed, under our assumptions, this problem is
solvable. On the other hand, given any solution $(x,v)$ of the system (\ref{xv.crit}), if the function $h$ is subdifferentially regular at the point $c(x)$ (as holds in particular if $h$ is convex or continuously differentiable), then we have the inclusion
\[
\nabla c(x)^* \partial h \big( c(x) \big) ~\subset~ \partial ( h \circ c)(x).
\]
Thus $0 \in \partial ( h \circ c)(x)$ and therefore $x$ must be a
critical point of the composite function $h \circ c$.

We can rewrite the criticality system \eqnok{xv.crit} in terms of the graph $\gph(\partial h)$ as follows:
\[
(c(x),v) \in \gph(\partial h), ~~ \nabla c(x)^* v = 0.
\]
Solving this system is often difficult in part because the graph $\gph(\partial h)$ may have a complicated structure.  Active set methods from classical nonlinear programming and its extensions essentially restrict attention to a suitable subset of $\gph(\partial h)$, thereby narrowing a local algorithmic search for a critical point.  We therefore make the following definition.
\begin{definition}
An {\bf actively sufficient set} for a critical point $\bar x$ of the composite function 
$h \circ c$ is a set $G \subset \gph(\partial h)$ containing a point of the form
$(c(\bar x),\bar v)$, where $\bar v$ is a multiplier vector for $\bar x$.
\end{definition}

The central idea we explore in this work is how to ``identify''  actively sufficient sets from among the parts of a decomposition of the graph $\gph(\partial h)$.  We present conditions ensuring that any sufficiently accurate approximate solution of system \eqnok{xv.crit} with the pair $\big[x, h\big(c(x)\big)\big]$ sufficiently near the pair $\big[\bar x, h\big(c(\bar x)\big)\big]$ identifies an actively sufficient set.  

\section{Main result}
We start with a useful tool.

\begin{lemma} \label{lem:vbd} 
Under Assumption \ref{blanket}, the set of multiplier vectors for $\bar x$ is nonempty and compact.
\end{lemma}
\begin{proof}
We have already observed the existence of a multiplier vector.  Since the subdifferential $\partial h \big( c(\bar x) \big)$ is a closed set, the set of multipliers must also be closed.  Assuming for contradiction that this set is unbounded, we can find
  a sequence $\{ v_r \}$ with $|v_r| \to \infty$ and 
\[
v_r \in \partial h \big( c(\bar x) \big), ~~ \nabla c(\bar x)^* \bar v_r = 0.
\]
By defining $w_r := v_r / |v_r|$, we have $|w_r| \equiv 1$ and hence
without loss of generality we can assume $w_r \to \bar{w}$ with
$|\bar{w}|=1$. Clearly, since $w_r \in N(\nabla c(\bar x)^*)$ and
the null space is closed, we have 
$\bar{w} \in N(\nabla c(\bar x)^*)$.  On the other hand,
$\bar{w} \in \partial^{\infty} h \big( c(\bar x) \big)$
follows from the definition of the horizon subdifferential. Since
$\bar{w} \neq 0$ we have a contradiction to condition \eqnok{cq}.
\end{proof}

We are ready to present the main result.

\begin{theorem} \label{MAIN}
Suppose Assumption \ref{blanket} holds.
Consider any closed set $G \subset \gph(\partial h)$
Then for any sufficiently small number $\epsilon > 0$, there exists a
number $\delta > 0$ with the following property.  For any point $x \in \R^n$ close to 
$\bar x$, in the sense that
\[
|x-\bar x| < \delta,
\]
if there exists a pair $(\hat{c},v) \in \R^m \times \R^m$ close to $G$, in the sense that
\[
\dist \big( (\hat{c},v),G \big) < \epsilon,
\]
and such that the first-order conditions hold approximately, in the sense that
\[
v \in \partial h (\hat{c}),
~~|\hat{c}-c(x)| < \delta,  
~~\big| h (\hat{c}) -  h \big( c(\bar x) \big) \big| < \delta,
~~\mbox{and}~~ | \nabla c(x)^* v | < \delta,
\]
then $G$ is an actively sufficient set for $\bar x$.
\end{theorem}

\begin{proof}
Suppose the result fails.  Then $G$ is not an actively sufficient set, and yet there exists a sequence of strictly
positive numbers $\epsilon_j \downarrow 0$ as $j \to \infty$ such
that, for each $j=1,2,\ldots$, the following property holds:  There
exist sequences
\[
x_j^r \in \R^n, ~~~ 
c_j^r \in \R^m, ~~~
v_j^r \in \partial h \big( c_j^r \big), ~~~ r=1,2,\ldots,
\]
satisfying
\[
x_j^r \to \bar x, ~~~ 
c_j^r \to c(\bar{x}), ~~~
h(c_j^r) \to h \big( c(\bar x) \big), ~~~ 
\nabla c(x_j^r)^* v_j^r \to 0,
\]
as $r \to \infty$, and yet
\[
\dist \big( (c_j^r,v_j^r),G \big) < \epsilon_j, ~~~ r=1,2,\ldots.
\]
For each $j$, we can use the proof
technique of Lemma~\ref{lem:vbd} to show that the sequence
$(v_j^r)_{r=1}^{\infty}$ must be bounded. Thus, by taking a
subsequence of the indices $r$, we can suppose that this sequence converges to
some vector $v_j$, which must be a multiplier vector at $\bar x$.
By continuity, we deduce
\[
\dist \big( (c(\bar x),v_j),G \big) \le \epsilon_j.
\]

By Lemma~\ref{lem:vbd}, the sequence $(v_j)_{j=1}^{\infty}$ is bounded, so after taking a subsequence of the indices
$j$, we can suppose that it converges to some multiplier vector $\bar v$.  Noting that the set $G$ is closed, we have by taking
limits as $j \to \infty$ that $(c(\bar x),\bar v) \in G$, contradicting the assumption that $G$ is not an actively sufficient set.
\end{proof}

An easy corollary extends from one potential actively sufficient set to many.

\begin{corollary} \label{maincor}
Suppose Assumption \ref{blanket} holds.  Consider any finite family ${\cal G}$ of closed subsets of $\gph(\partial h)$.  Then for any sufficiently small number $\epsilon > 0$, there exists a number $\delta > 0$ with the following property.  For any point $x \in \R^n$ close to $\bar x$, in the sense that
\begin{equation} \label{close}
|x-\bar x| < \delta,
\end{equation}
if there exists a pair $(\hat{c},v) \in \R^m \times \R^m$ close to some set $G \in {\cal G}$, in the sense that
\begin{equation} \label{sufficient}
\dist \big( (\hat{c},v),G \big) < \epsilon,
\end{equation}
such that the first-order conditions hold approximately, in the sense that
\begin{equation} \label{approximate}
v \in \partial h (\hat{c}),
~~|\hat{c}-c(x)| < \delta,  
~~\big| h (\hat{c}) -  h \big( c(\bar x) \big) \big| < \delta,
~~\mbox{and}~~ | \nabla c(x)^* v | < \delta,
\end{equation}
then $G$ is an actively sufficient set for $\bar x$.
\end{corollary}

\begin{proof}
For each set $G \in {\cal G}$, we apply Theorem \ref{MAIN}, deducing the existence of a number $\epsilon_G > 0$ such that the conclusion of the theorem holds for all numbers $\epsilon$ in the interval $(0,\epsilon_G)$.  Define the strictly positive number $\bar \epsilon = \min_G \epsilon_G$.  We claim the result we seek holds for all $\epsilon$ in the interval 
$(0,\bar \epsilon)$.  To see this, we apply the theorem for each set $G \in {\cal G}$ to deduce the existence of a number 
$\delta_G > 0$ such that the conditions \eqnok{approximate} and \eqnok{close}, with 
$\delta = \delta_G$, and the condition \eqnok{sufficient}, together imply that $G$ is a actively sufficient set for $\bar x$.  The result now follows by setting 
$\delta = \min_G \delta_G$.
\end{proof}

The following result is a simple special case, easily proved directly.

\begin{corollary}
Under the assumptions of Corollary \ref{maincor},
there exists a number
$\tilde{\epsilon}>0$ such that
\begin{equation} \label{galah}
\dist \big( (c(\bar{x}),\bar{v}), G \big) > \tilde{\epsilon}
\end{equation}
for all multiplier vectors $\bar v$ for the critical point $\bar x$, and all sets 
$G \in {\cal G} $ that are not actively sufficient for $\bar x$.
\end{corollary}
\begin{proof}
In Corollary \ref{maincor}, set $x=\bar x$ and $\hat c = c(\bar x)$.
\end{proof}

We end this section with another corollary, indicating how we might use the main result in practice.

\begin{corollary} \label{sequence}
Suppose Assumption \ref{blanket} holds.  Consider any finite family ${\cal G}$ of closed subsets of $\gph(\partial h)$.  Then for any sequence of points $x_r \in \R^n$, vectors $c_r \in \R^m$, subgradients $v_r \in \partial h(c_r)$, and sets $G_r \in {\cal G}$ (for $r=1,2,\ldots$), satisfying
\begin{eqnarray*}
&
x_r \to \bar x, ~~~ 
|c_r - c(x_r)| \to 0,~~~ 
h(c_r) \to h\big(c(\bar x)\big),
& \\
&
\nabla c(x_r)^* v_r \to 0,~~~
\dist\big( (c_r,v_r),G_r \big) \to 0,
&
\end{eqnarray*}
as $r \to \infty$, the set $G_r$ is actively sufficient for $\bar x$
for all $r$ sufficiently large.
\end{corollary}
\begin{proof}
  Apply Corollary \ref{maincor}, for any sufficiently small number
  $\epsilon > 0$.  Then, for the number $\delta > 0$ guaranteed by the
  corollary, equations \eqnok{close}, \eqnok{sufficient} and
  \eqnok{approximate} hold for all sufficiently large $r$, so the
  result follows.
\end{proof}

\section{Subdifferential graph decomposition}
To apply the ideas in the previous section, we typically assume the
availability of a decomposition of $\gph (\partial h)$ (the graph of
the subdifferential of $h$) into some finite union of closed, not
necessarily disjoint sets $G^1,G^2,\ldots,G^k \subset \R^m \times
\R^m$.  For this decomposition to be useful, the sets $G^i$ should be
rather simple, so that the restricted system
\[
(c(x),v) \in G^i, ~~ \nabla c(x)^* v = 0.
\]
is substantially easier to solve than the original criticality system.
The more refined the
decomposition, the more information we may be able to derive from the
identification process.  Often we have in mind the situation where
each of the sets $G^i$ is a polyhedron. We might, for example, assume
that whenever some polyhedron is contained in the list $(G^i)$, so is
its entire associated lattice of closed faces.

\begin{Exa}[Scalar examples] \label{ex:scalarh}
{\rm
We give some simple examples in the case $m=1$.  Consider first the
indicator function for $\R_+$, defined by $h(c)=0$ for $c\ge 0$ and
$+\infty$ for $c<0$. We have 
\[
\partial h(c) =  \begin{cases}
\emptyset & \mbox{\rm if $c<0$} \\
(-\infty,0] & \mbox{\rm if $c=0$} \\
\{0\} & \mbox{\rm if $c>0$.}
\end{cases}
\]
Thus an appropriate decomposition is $\gph (\partial h) = G^1 \cup
G^2 \cup G^3$, where
\[
G^1 = \{0 \} \times (-\infty,0], \qquad
G^2 = \{(0,0)\}, \qquad
G^3 = [0,\infty) \times \{0\}.
\]
Similar examples are the absolute value function $|\cdot|$, for which a
decomposition is $\gph (\partial \, |\cdot|) = G^1 \cup G^2 \cup G^3$,
where
\begin{equation} \label{abs}
G^1 = (-\infty,0] \times \{-1\}, \qquad
G^2 = \{0 \} \times [-1,1], \qquad
G^3 = [0,\infty) \times \{1\}
\end{equation}
(further refinable by including the two sets $\{0,\pm 1\}$),
and the positive-part function $\mbox{\rm pos}(c) = \max(c,0)$, for which a
decomposition is $\gph (\partial \, \mbox{\rm pos}) = G^4 \cup G^5 \cup G^6$,
where
\begin{equation} \label{pos}
G^4 = (-\infty,0] \times \{0\}, \qquad
G^5 = \{0 \} \times [0,1], \qquad
G^6 = [0,\infty) \times \{1\}
\end{equation}
(again refinable).
A last scalar example, which involves a nonconvex function $h$, is given by
$h(c) = 1- e^{-\alpha |c|}$ for some constant $\alpha>0$. We have
\[
\partial h(c) = \begin{cases}
\{-\alpha e^{\alpha c}\} & \mbox{\rm if $c<0$} \\
[-\alpha,\alpha] & \mbox{\rm if $c=0$} \\
\{\alpha e^{-\alpha c}\} & \mbox{\rm if $c>0$}.
\end{cases}
\]
An appropriate partition is $\gph (\partial h) = G^1 \cup G^2 \cup G^3$,
where
\[
G^1 = \{(c,-\alpha e^{\alpha c}) : c \le 0 \} \qquad
G^2 = \{0\} \times [-\alpha,\alpha] \qquad
G^3 = \{(c,\alpha e^{-\alpha c}) : c \ge 0 \}.
\]
}
\end{Exa}

\begin{Exa}[An $\ell_1$-penalty function] \label{ex:2dh} 
{\rm
Consider a function $h:\R^2 \to \R$ that is an
  $\ell_1$-penalty function for the constraint system $c_1=0$, $c_2
  \le 0$, that is,
  \begin{equation}
    \label{eq:h2l1}
    h(c) = |c_1| + \max(c_2,0).
  \end{equation}
Using the notation of the previous example, we have 
\[
\partial h(c_1,c_2) = \partial(|\cdot|)(c_1) \times \partial\mbox{\rm pos}(c_2).
\]
A partition of $\gph (\partial h)$ into nine closed
sets  can be constructed by using interleaved Cartesian products of
\eqnok{abs} and \eqnok{pos}.
}
\end{Exa}

Much interest lies in the
case in which the function $h$ is polyhedral, so that 
$\gph (\partial h)$ is a finite union of polyhedra.  However, the latter property
holds more generally for the ``piecewise linear-quadratic'' functions
defined in \cite{Roc98}.  

Of course, we cannot decompose the graph of the subdifferential
$\partial h$ into a finite union of closed sets unless this graph is
itself closed.  This property may fail, even for quite simple
functions.  For example, the lower semicontinuous function $h \colon
\R \to \R$ defined by $h(c) = 0$ for $c \le 0$ and $h(c) = 1-c$ for
$c>0$ has subdifferential given by
\[
\partial h(c) =  \begin{cases}
\{0\} & \mbox{\rm if $c<0$} \\
[0,\infty) & \mbox{\rm if $c=0$} \\
\{-1\} & \mbox{\rm if $c>0$,}
\end{cases}
\]
so $\gph (\partial h)$ is not closed.  On the other hand, the subdifferentials of lower semicontinuous {\em convex} functions are closed.

In general, for any semi-algebraic
function $h$, the set $\gph(\partial h)$ is semi-algebraic.  If this set is also closed, then it stratifies into a finite union of smooth manifolds with boundaries.
In concrete cases, a decomposition may be reasonably straightforward.  We end this section with two examples.

\begin{Exa}
{\rm
The graph of the subdifferential of the Euclidean norm on $\R^n$ decomposes into the union of the following two closed sets:
\[
\{(0,v) : |v| \le 1\} ~~~\mbox{and}~~~ 
\Big\{ \Big(c,\frac{1}{|c|}c \Big) : c \neq 0 \Big\} ~\cup~ 
\big\{(0,v) : |v| = 1 \big\}.
\]
}
\end{Exa}

\begin{Exa}
{\rm
Consider the maximum eigenvalue function $\lambda_{\mbox{\scriptsize max}}$
on the Euclidean space ${\mathbf S}^k$ of $k$-by-$k$ symmetric matrices (with the inner product $\ip{X}{Y} = \mbox{trace}(XY)$).  In this space, the following sets are closed:
\begin{eqnarray*}
{\mathbf S}^k_r & = & \{ Y \in {\mathbf S}^k : Y~ \mbox{has rank} \le r \}
~~~ (r = 0,1,\ldots,k)  \\
\mbox{}_m{\mathbf S}^k
 & = & \{ X \in {\mathbf S}^k : 
\lambda_{\mbox{\scriptsize max}}(X)~ \mbox{has multiplicity} \ge m \}
~~~ (m = 1,2,\ldots,k).
\end{eqnarray*}
Trivially we can decompose the graph
$\mbox{gph}(\partial \lambda_{\mbox{\scriptsize max}})$ into its intersection with each of the sets $\mbox{}_m{\mathbf S}^k \times {\mathbf S}^k_r$.  However, we can simplify, since it is is well known (see \cite{Lew96Hermitian}, for example) that $\partial \lambda_{\mbox{\scriptsize max}}(X)$ consists of matrices of rank no more than the multiplicity of $\lambda_{\mbox{\scriptsize max}}(X)$.  Hence we can decompose the graph into the union of the sets
\[
G_{m,r} ~=~
\mbox{gph}(\partial \lambda_{\mbox{\scriptsize max}}) 
\cap 
(\mbox{}_m{\mathbf S}^k \times {\mathbf S}^k_r)
~~~ (1 \le r \le m \le k).
\]
To apply the theory we have developed, we need to measure the distance from any given pair $(X,Y)$ in the graph to each of the sets $G_{m,r}$.  This is straightforward, as follows.  A standard characterization of $\partial \lambda_{\mbox{\scriptsize max}}$ \cite{Lew96Hermitian} shows that there must exist an orthogonal matrix $U$, a vector $x \in \R^k$ with nonincreasing components, and a vector $y \in \R^k_+$ satisfying $\sum_i y_i = 1$ and $y_i = 0$ for all indices $i > p$, where $p$ is the multiplicity of the largest component of $x$, such that the following simultaneous spectral decomposition holds:  $X = U^T(\mbox{Diag}\,x) U$ and $Y = U^T (\mbox{Diag}\,y) U$.  Now define a vector $\tilde x \in \R^k$ by replacing the first $m$ components of $x$ by their mean.  (Notice that the components of $\tilde x$ are then still in nonincreasing order, and the largest component has multiplicity at least $p$.)  Define a vector $\tilde y \in \R^k$ by setting all but the largest $r$ components of $y$ to zero and then rescaling the resulting vector to ensure its components sum to one.  (Notice that $\tilde y_i = 0$ for all indices 
$i > p$.)  Finally, define matrices $\tilde X = U^T(\mbox{Diag}\,\tilde x) U$ and 
$\tilde Y = U^T (\mbox{Diag}\,\tilde y) U$.  Then, by the same subdifferential characterization, we have $\tilde Y \in \partial \lambda_{\mbox{\scriptsize max}}(\tilde X)$, so in fact $(\tilde X, \tilde Y) \in G_{m,r}$.  Hence the distance from $(X,Y)$ to $G_{m,r}$ is at most $\sqrt{|x-\tilde x|^2 + |y-\tilde y|^2}$.  In fact this easily computable estimate is exact, since it is well known that $\tilde Y$ is a closest matrix to $Y$ in the set ${\mathbf S}^k_r$ and, by \cite[Example A.4]{Lew08aAlternating}, $\tilde X$ is a closest matrix to $X$ in the set $\mbox{}_m{\mathbf S}^k$. 
} 
\end{Exa}

\section{Classical nonlinear programming} \label{sec:classical}
We illustrate all of our key concepts on the special case of classical nonlinear programming, which we state as follows:
\[
\mbox{\rm (NLP)} \hspace{2cm}
\left\{
\begin{array}{lrcll}
\inf              &   f(x) &   &       &                  \\
\mbox{subject to} & p_i(x) & = & 0     & (i=1,2,\ldots,s) \\
                  & q_j(x) &\le& 0     & (j=1,2,\ldots,t) \\
                  &     x  &\in& \R^n, &
\end{array}
\right.
\]
where the functions $f,p_i,q_j \colon \R^n \to \R$ are all
continuously differentiable.  We use the notation
\begin{equation} \label{qpm}
q^+(x) = \max(q(x),0), ~~~~ q^-(x) = \min(q(x),0),
\end{equation}
where the max and min of $q(x) \in \R^t$ are taken componentwise. (It
follows that $q(x) = q^+(x)+q^-(x)$.)

We can model the problem (NLP) in our
composite form (\ref{hc}) by defining a continuously differentiable
function $c \colon \R^n \to \R \times \R^s \times \R^t$ and a
polyhedral function $h \colon \R \times \R^s \times \R^t \to \bar\R$
through
\begin{subequations} \label{nlp}
\begin{align}
\label{nlp.1}
c(x) & = \big( f(x),p(x),q(x) \big) ~~~ (x \in \R^n) \\
\label{nlp.2}
h(u,y,w) & =  
\left\{
\begin{array}{ll}
u & (y=0,~ w \le 0) \\
+\infty & (\mbox{otherwise})
\end{array}
\right.
~~~ (u \in \R,~ y \in \R^s,~ w \in \R^t).
\end{align}
\end{subequations}
Clearly for any point $x \in \R^n$, the adjoint map $\nabla c(x)^* \colon \R \times \R^s \times \R^t \to \R^n$ is given by
\[
\nabla c(x)^*(\theta,\lambda,\mu) 
~=~ \theta \nabla f(x) + \sum_i \lambda_i \nabla p_i(x) + \sum_j \mu_j \nabla q_j(x).
\]
The subdifferential and horizon subdifferential of $h$ at any point 
$(u,0,w) \in \R \times \R^s \times \R_-^t$ are given by
\begin{eqnarray*}
\partial h(u,0,w) & = & \{1\} \times \R^s \times \{ \mu \in \R_+^t : \ip{\mu}{w} = 0 \} \\
\partial^{\infty} h(u,0,w) & = & 
\{0\} \times \R^s \times \{ \mu \in \R_+^t : \ip{\mu}{w} = 0 \}.
\end{eqnarray*}
(Elsewhere in $\R \times \R^s \times \R^t$, these two sets are
respectively $\emptyset$ and $\{0\}$.)

Armed with these calculations, consider any critical point $\bar x$
(or in particular, any local minimizer for the nonlinear program).  By
assumption, $\bar x$ is a feasible solution.  Classically, the {\em active
set} is
\[
\bar J ~=~ \{j : q_j(\bar x) = 0 \}.
\]
The regularity condition \eqref{cq} becomes the following assumption.
\begin{assumption}[Mangasarian-Fromovitz] \label{mf}
The only pair $(\lambda,\mu) \in \R^s \times \R_+^t$ satisfying $\mu_j = 0$ for 
$j \not\in \bar J$ and
\[
\sum_i \lambda_i \nabla p_i(\bar x) + \sum_j \mu_j \nabla q_j(\bar x) = 0
\]
is $(\lambda,\mu) = (0,0)$.
\end{assumption}

\noindent

In this framework, what we have called a multiplier vector for the
critical point $\bar x$ is just a pair $(\bar\lambda,\bar\mu) \in \R^s
\times \R^t_+$ satisfying $\bar\mu_j = 0$ for $j \not\in \bar J$ and
\begin{equation} \label{stationary}
\nabla f(\bar x) + \sum_i \bar\lambda_i \nabla p_i(\bar x) + 
\sum_j \bar\mu_j \nabla q_j(\bar x) = 0.
\end{equation}
It is evident that Lemma~\ref{lem:vbd} retrieves the classical
first-order optimality conditions: existence of Lagrange multipliers
under the Mangasarian-Fromovitz constraint qualification.

Nonlinear programming is substantially more difficult than solving
nonlinear systems of equations, because we do not know the active set $\bar J$ in advance.  Active set methods try to identify $\bar J$,
since, once this set is know, we can find a stationary point by
solving the system
\begin{eqnarray*}
\nabla f(x) + \sum_i \lambda_i \nabla p_i(x) + \sum_{j \in \bar J} \mu_j \nabla q_j(x) 
& = & 0 \\
p_i(x) & = & 0 ~~~ (i=1,2,\ldots,p) \\
q_j(x) & = & 0 ~~~ (j \in \bar J),
\end{eqnarray*}
which is a nonlinear system of $n + p + |\bar J|$ equations for the vector
$(x,\lambda,\mu_J) \in \R^n \times \R^p \times \R^{|J|}$.  
Our aim here is to formalize this process of identification.  
Our approach broadly follows that of \cite{ObeW05a}, with extensive
generalization to the broader framework of composite minimization.

The classical notion of active set in nonlinear programming
arises from a certain combinatorial structure in the graph of the
subdifferential $\partial h$ of the outer function $h$:
\begin{equation} \label{gphh}
\gph(\partial h) ~=~ 
\big\{ \big( (u,0,w) , (1,\lambda,\mu) \big) : w \le 0,~ \mu \ge 0,~ \ip{w}{\mu} = 0 \big\}.
\end{equation}
We can decompose this set into a finite union of polyhedra, as follows:
\[
\gph(\partial h) = \bigcup_{J \subset \{1,2,\ldots,t\}} G^J,
\]
where
\begin{equation} \label{def.GJ}
G^J 
=
\Big\{ \big( (u,0,w) , (1,\lambda,\mu) \big) : 
w \le 0,~ \mu \ge 0,~ w_j=0~ (j \in J),~ \mu_j=0~ (j \not\in J) \Big\}.
\end{equation}
According to our definition, $G^J$ is an actively sufficient set exactly when 
$J \subset \bar J$ and there exist vectors $\bar\lambda \in \R^s$ and 
$\bar\mu \in \R^t_+$ satisfying $\bar\mu_j = 0$ for all $j \not\in J$, and the stationarity condition \eqnok{stationary}.  We call such an index set $J$ {\em sufficient} at $\bar x$.

We next illustrate the main result.  We use the notation \eqnok{qpm}
below. In addition, for a vector $q \in \R^t$ and a nonnegative scalar
$\delta$, we define $q^{\delta} \in \R^t$ as follows:
\begin{equation} \label{qdel}
q^{\delta}_i = \begin{cases} q_i, & \mbox{if} ~ q_i < -\delta, \\
0, & \mbox{if} ~ q_i \ge -\delta.
\end{cases}
\end{equation}

\begin{corollary}
Consider a critical point $\bar x \in \R^n$ for the nonlinear program
{\rm (NLP)}, where the objective function and each of the constraints
functions are all continuously differentiable.  Suppose the
Mangasarian-Fromovitz condition (Assumption \ref{mf}) holds.  Then for
any sufficiently small number $\epsilon' > 0$, there exists a number
$\delta' > 0$ with the following property.  For any triple
$(x,\lambda,\mu) \in \R^n \times \R^s \times \R^t_+$ satisfying 
\begin{subequations} \label{dog}
\begin{align}
\label{dog.1}
|x - \bar x| & < \delta', \\
\label{dog.2}
\mu_j = 0 ~\mbox{whenever}~ q_j(x) & < -\delta', \\
\label{dog.3}
\Big| \nabla f(x) + \sum_{i=1}^s \lambda_i \nabla p_i(x) + \sum_{j=1}^t \mu_j \nabla q_j(x) \Big|
 & < \delta',
\end{align}
\end{subequations}
any index set $J \subset \{1,2,\ldots,t\}$ that satisfies
\begin{subequations} \label{approx-comp}
\begin{alignat}{2} \label{approx-comp.1}
q_j(x) & > -\epsilon' & \qquad & \mbox{\rm for all $j \in J$}, \\
\label{approx-comp.2}
\mu_j & < \epsilon'  & \qquad & \mbox{\rm for all $j \not\in J$},
\end{alignat}
\end{subequations}
is sufficient for $\bar x$.
\end{corollary}

\begin{proof}
Applying Corollary~\ref{maincor} using the decomposition above, 
for any number $\epsilon > 0$
sufficiently small, there exists a number $\delta > 0$ with the
following property. For any $(x, \theta, \lambda, \mu, \hat{f},
\hat{p}, \hat{q}) \in \R^n \times \R \times \R^s \times \R^t \times \R
\times \R^s \times \R^t$ such that
\begin{subequations} \label{poodle}
\begin{align}
\label{poodle.1}
|x-\bar{x}| & ~<~ \delta, \\
\label{poodle.2}
\big| (\hat{f},\hat{p},\hat{q}) - \big(f(x),p(x),q(x)\big) \big| &~<~ \delta, \\
\label{poodle.3}
\big| h(\hat{f},\hat{p},\hat{q}) - h\big(f(\bar{x}),p(\bar{x}),q(\bar{x})\big) \big| &~<~ \delta, \\
\label{poodle.4}
\Big| \theta \nabla f(x) + \sum_{i=1}^s \lambda_i \nabla p_i(x) + \sum_{j=1}^t \mu_j \nabla q_j(x) \Big| &~<~ \delta, \\
\label{poodle.5}
(\theta,\lambda,\mu) & ~\in~ \partial h(\hat{f},\hat{p},\hat{q}),
\end{align}
\end{subequations}
and for any index set $J \subset \{1,2,\dotsc,t \}$ such that
\begin{equation} \label{sheep}
\dist \Big( \big( (\hat{f},\hat{p},\hat{q}), (\theta,\lambda,\mu) \big),
G^J \Big) < \epsilon,
\end{equation}
we have that there exist multipliers $\bar{\lambda} \in \R^s$ and
$\bar{\mu} \in \R^t_+$ such that
\begin{subequations} \label{schnauzer}
\begin{align}
\label{schnauzer.1}
\left((f(\bar{x}),p(\bar{x}),q(\bar{x})), (1,
\bar{\lambda},\bar{\mu}) \right) & \in G^J, \\
\label{schnauzer.2}
\nabla f(\bar x) + \sum_{i=1}^s \bar\lambda_i \nabla p_i(\bar x) + 
\sum_{j=1}^t \bar\mu_j \nabla q_j(\bar x) &= 0.
\end{align}
\end{subequations}

To prove our claim, we need to perform three tasks.
\begin{itemize}
\item[(i)] Identify a value of $\delta'$ and values of $\theta$ and
  $(\hat{f},\hat{p},\hat{q})$ such that \eqnok{poodle} holds whenever
  $(x,\lambda,\mu)$ satisfies \eqnok{dog};
\item[(ii)] Identify a value of $\epsilon'$ such that for these
  choices of $x$, $(\hat{f},\hat{p},\hat{q})$, and
  $(\theta,\lambda,\mu)$, the condition \eqnok{approx-comp} implies
  that \eqnok{sheep} is satisfied.
\item[(iii)] Prove that the outcome of Corollary~\ref{maincor}, namely
  \eqnok{schnauzer}, implies that the index set $J$ is sufficient.
\end{itemize}

We start with (i). We choose $\delta'>0$ to satisfy $\delta' \le
\delta$ and $\delta' \le \epsilon/\sqrt{t}$, and also small enough
that $|x-\bar{x}| < \delta'$ implies
\begin{subequations} \label{pigdog}
\begin{align}
\label{pigdog.1}
|p(x)-p(\bar{x})| + \sqrt{t} \delta' + |q^+(x)-q^+(\bar{x})| & < \delta, \\
\label{pigdog.2}
|f(x)-f(\bar{x})| & < \delta.
\end{align}
\end{subequations}
Now set $\theta=1$ and $(\hat{f},\hat{p},\hat{q}) =
(f(x),0,q^{\delta'}(x))$. Note that by \eqnok{dog.2} and \eqnok{qdel},
$\mu_j=0$ whenever $q_j^{\delta'}(x) \neq 0$, and $\mu_j \ge 0$
otherwise.  We thus have from \eqnok{gphh} that $(1,\lambda,\mu) \in
\partial h (f(x),0,q^{\delta'}(x))$, so that \eqnok{poodle.5}
holds. Since $\delta' \le \delta$ and $|x-\bar{x}| < \delta'$, we have
\eqnok{poodle.1} immediately, while \eqnok{poodle.4} follows from
$\theta=1$ and \eqnok{dog.3}.

We have from $p(\bar{x})=0$ and $q^+(\bar{x})=0$ that 
\begin{align*}
\big| (\hat{f},\hat{p},\hat{q}) & - (f(x),p(x),q(x)) \big| \\
& \le |p(x)| + \left| q^{\delta'}(x)-q(x) \right| \\
& \le |p(x)-p(\bar{x})| + \left| q^{\delta'}(x)-q^-(x) \right| +
\left| q^+(x)-q^+(\bar{x}) \right| \\
& \le |p(x)-p(\bar{x})| + \sqrt{t} \delta' + 
\left| q^+(x)-q^+(\bar{x}) \right|  < \delta,
\end{align*}
by \eqnok{pigdog.1}, so that \eqnok{poodle.2} holds. Further, 
\begin{align*}
\big| h(\hat{f},\hat{p},\hat{q}) & - h(f(\bar{x}),p(\bar{x}),q(\bar{x})) \big| \\
&= \left| h(f(x),0,q^{\delta'}(x)) - h(f(\bar{x}),p(\bar{x}),q(\bar{x})) \right| \\
&= |f(x)-f(\bar{x})| < \delta,
\end{align*}
by \eqnok{pigdog.2}, so that \eqnok{poodle.3} holds. At this point we
have completed task (i).

We now show (ii). Define $\epsilon' = \epsilon / \sqrt{t}$ and note
that by one of our conditions on $\delta'$, we have $\delta' \le
\epsilon'$.  Defining vectors $\hat w,\hat \mu \in \R^t$ by
\begin{equation} \label{porkers}
\hat w_j = 
\begin{cases} 
q_j(x) & \mbox{if}~ q_j(x) \le -\epsilon' \\
0      & \mbox{otherwise,}
\end{cases}
\hspace{15mm}
\hat \mu_j = 
\begin{cases}
\mu_j & \mbox{if}~ \mu_j \ge \epsilon' \\
0      & \mbox{otherwise,}
\end{cases}
\end{equation}
then by \eqnok{approx-comp} and  the definition of $G^J$ we have
\[
\left( ( f(x),0,\hat w ) , (1,\lambda,\hat \mu) \right) \in G^J.
\]
Since $\delta' \le \epsilon'$, we have
\begin{equation} \label{goats}
q^{\delta'}_i(x)-\hat{w}_i = 
\begin{cases}
0, & \mbox{if} ~~ q_i(x) \le -\epsilon', \\
q_i(x), & \mbox{if} ~~ q_i(x) \in (-\epsilon',-\delta'), \\
0, & \mbox{if} ~~ q_i(x) \ge -\delta'.
\end{cases}
\end{equation}
Thus in \eqnok{sheep}, using the values of $(\hat{f},\hat{p},\hat{q})$
and $(\theta,\lambda,\mu)$ defined above, we have that
\begin{align*}
\dist \Big( \Big( & \big(f(x),0, q^{\delta'}(x) \big), (1,\lambda,\mu) \Big),G^J \Big)^2 \\
&\le
\left| 
\left( \big( f(x),0,q^{\delta'}(x) \big) , (1,\lambda, \mu) \right) -
\left( ( f(x),0,\hat w ), (1,\lambda,\hat \mu) \right)
\right|^2
\\
&  =
\left|q^{\delta'}(x) - \hat w \right|^2 + |\mu - \hat \mu|^2 ~\le~ t (\epsilon')^2 = \epsilon^2,
\end{align*}
The final inequality in this expression follows from (\ref{porkers})
and (\ref{goats}) together with the fact that we cannot have both
$q^{\delta'}_i(x) \neq \hat{w}_i$ and $\mu_i \neq \hat{\mu}_i$ for any
index $i$. If $q^{\delta'}_i(x) \neq \hat{w}_i$, we have from
\eqnok{goats} that $q^{\delta'}_i(x) \in (-\epsilon',-\delta')$, thus
$\mu_i=0$ by \eqnok{dog.2}, thus $\hat{\mu}_i=0$ by \eqnok{porkers},
thus $|\mu_i - \hat{\mu}_i|=0$.  We conclude that the inequality
\eqnok{sheep} is satisfied, completing the proof of part (ii).

Part (iii) of the proof is immediate from the definition of a
sufficient index set, so the proof is complete.
\end{proof}

\section{Partial smoothness}
We next observe a connection between the decomposition ideas we have
introduced and the notion of ``partial smoothness''
\cite{Lew03Active}.  For simplicity, in this section we restrict to
the convex case, although extensions are possible.  A lower
semicontinuous convex function $h \colon \R^m \to \bar\R$ is {\em
  partly smooth} at point $\bar c \in \R^m$ relative to a set ${\cal
  M}$ containing $\bar c$ when ${\cal M}$ is a manifold around $\bar
c$, the restricted function $h|_{\cal M}$ is $C^2$, and the
subdifferential mapping $\partial h$ is continuous at $\bar c$ when
restricted to ${\cal M}$ with $\partial h(\bar c)$ having affine span
a translate of the normal space to ${\cal M}$ at $\bar c$.

\begin{theorem} \label{partly} Consider a lower semicontinuous convex
  function $h \colon \R^m \to \bar\R$, a point $\bar c \in \R^m$, and
  a vector $\bar v$ lying in the relative interior of the
  subdifferential $\partial h(\bar c)$.  Suppose that $h$ is partly
  smooth at $\bar c$ relative to a closed set ${\cal M} \subset \R^m$.
  Then the graph of the subdifferential $\partial h$ is the union of
  the following two closed sets:
\[
G^1 = \{ (c,v) : c \in {\cal M},~ v \in \partial h(c) \}, ~~~
G^2 = \mbox{\rm cl} \{ (c,v) : c \not\in {\cal M},~ v \in \partial h(c) \}.
\]
Furthermore, the set $G^2$ does not contain the point $(\bar c, \bar v)$.
\end{theorem}
\begin{proof}
  As is well known, since $h$ is convex and lower semicontinuous,
  $\mbox{gph}(\partial h)$ is closed: indeed we can write it as the
  lower level set of a lower semicontinuous function:
\[
\mbox{gph}(\partial h) ~=~ \{ (c,v) : h(c) + h^*(v) - \ip{c}{v} \le 0 \},
\]
where $h^*$ denotes the Fenchel conjugate of $h$.  Since the set $G^1$
is just $\mbox{gph}(\partial h) \cap ({\cal M} \times \R^m)$, and
since ${\cal M}$ is closed by assumption, $G^1$ is a closed subset of
the graph $\mbox{gph}(\partial h)$.  The set $G^2$ is closed by
definition, and $G^2$ is also obviously a subset of $\gph (\partial
h)$.  Therefore, we have the decomposition $\mbox{gph}(\partial h) =
G^1 \cup G^2$.

It remains to show $(\bar c, \bar v) \not\in G^2$.  If this property
fails, then there is a sequence of points $c_r \not\in {\cal M}$
($r=1,2,\ldots$) approaching the points $\bar c$, and a corresponding
sequence of subgradients $v_r \in \partial h(c_r)$ approaching the
subgradient $\bar v$.  Then a standard subdifferential continuity
argument shows $h(c_r) \to h(\bar c)$: to be precise, we have
\begin{align*}
\liminf_r h(c_r) 
& = 
\liminf_r (\ip{c_r}{v_r} - h^*(v_r)) = \ip{\bar c}{\bar v} - \limsup_r h^*(v_r) \\
& \ge 
\ip{\bar c}{\bar v} - h^*(\bar v) = h(\bar c) \ge \limsup_r h(c_r).
\end{align*}
Now \cite[Thm 6.11]{LewW08a} implies the contradiction $c_r \in {\cal M}$ for all large $r$.
\end{proof}

We illustrate by showing how partial smoothness leads to identification.

\begin{corollary} \label{strict}
Suppose Assumption \ref{blanket} holds.
Suppose that the critical point $\bar x$ has a unique multiplier vector $\bar v$, and that 
$\bar v \in \mbox{\rm ri}\, \partial h \big( c(\bar x) \big)$.
Finally, assume that $h$ is convex, and partly smooth at the point $c(\bar x)$ relative to a closed set 
${\cal M} \subset \R^m$.  Then any sufficiently accurate solution of the criticality conditions near $\bar x$ must identify the set ${\cal M}$.  More precisely, for any sequence of points 
$x_r \in \R^n$, vectors $c_r \in \R^m$, and subgradients $v_r \in \partial h(c_r)$ (for $r=1,2,\ldots$), satisfying
\[
x_r \to \bar x, ~~~ 
|c_r - c(x_r)| \to 0,~~~ 
h(c_r) \to h\big(c(\bar x)\big),~~~
\nabla c(x_r)^* v_r \to 0,
\]
as $r \to \infty$, we must have $c_r \in {\cal M}$ for all sufficiently large $r$.
\end{corollary}
\begin{proof}
Consider the decomposition described in Theorem \ref{partly}.  Our assumptions imply that the set $G^2$ is not actively sufficient.  We now apply Corollary \ref{sequence} to deduce the result.
\end{proof}

\section{Identifying Activity via a Proximal Subproblem}
\label{sec:prox}

In this section we consider the question of whether closed sets $G$
that are actively sufficient at a solution $\bar{x}$ of the composite minimization problem
\eqnok{hc} can be identified from a nearby point $x$ by solving the
following subproblem:
\begin{equation} \label{pls} 
\min_{d} \, h_{x,\mu}(d) := 
h\big(c(x)+ \nabla c(x) d\big) + \frac{\mu}{2} |d|^2. 
\end{equation}
Properties of local solutions of this subproblem and of a first-order
algorithm based on it have been analyzed by the authors in
\cite{LewW08a}. In that work, we gave conditions guaranteeing in
particular that if the function $h$ is partly smooth relative to some
manifold ${\cal M}$ containing the critical point $\bar x$, then the
subproblem \eqnok{pls} ``identifies'' ${\cal M}$: that is, nearby
local minimizers must lie on ${\cal M}$.

The identification result from \cite{LewW08a} requires a rather strong
regularity condition at the critical point $\bar{x}$.  When applied to
the case of classical nonlinear programming we described above, this
condition reduces to the linear independence constraint qualification,
in particular always implying uniqueness of the multiplier vector.  In
the simplest case, when, in addition, strict complementarity holds,
there is a unique sufficient index set, in the terminology of Section
\ref{sec:classical}, and the identification result Corollary
\ref{strict} applies.

By contrast, in this section, we pursue more general identification
results, needing only the transversality condition \eqnok{cq}. Certain
additional assumptions on the function $h$ are required, whose purpose
is essentially to ensure that the solution of \eqnok{pls} is well
behaved.

We start with some technical results from \cite{LewW08a}, and then state our main result.

\begin{definition} \label{def:proxreg}
A function $h:\R^m \to \bar \R$ is {\bf prox-regular} at a point $\bar{c} \in \R^m$ if the value $h(\bar c)$ is finite and every point in $\R^m \times \R$ sufficiently close to the point $\big(\bar c, h(\bar c)\big)$ has a unique nearest point in the epigraph 
$\{(c,t) : t \ge h(c) \}$.
\end{definition}

\noindent
In particular, lower semicontinuous convex functions are everywhere prox-regular, as are sums of continuous convex functions and $\cC^2$ functions.

For the results that follow, we need to strengthen our underlying Assumption \ref{blanket}, as follows.

\begin{assumption} \label{blanket2}
In addition to Assumption \ref{blanket}, the function $c$ is $\cC^2$ around the critical point $\bar x$ and the function $h$ is prox-regular at the point $\bar c = c(\bar x)$.
\end{assumption}

The following result is a restatement of \cite[Theorem~6.5]{LewW08a}.  It concerns
existence of local solutions to \eqnok{pls} with nice properties.
\begin{theorem} \label{th:proxstep} 
Suppose Assumption \ref{blanket2} holds.  Then there exist
  numbers $\bar\mu \ge 0$, $\delta>0$ and $k>0$ and a mapping 
  $d \colon B_{\delta}(\bar{x}) \times (\bar{\mu},\infty) \to \R^n$ such that
  the following properties hold.
\begin{itemize}
\item[{\rm (a)}]
For all points $x \in B_{\delta}(\bar{x})$ and all scalars $\mu > \bar{\mu}$,
the point $d(x,\mu)$ is a local minimizer of the subproblem (\ref{pls}), and moreover satisfies $|d(x,\mu)| \le k |x-\bar{x}|$.

\item[{\rm (b)}] Given any sequences of points $x_r \to \bar{x}$ and scalars $\mu_r >
  \bar{\mu}$, if either $h(c(x_r)) \to h(\bar{c})$ or $\mu_r |x_r-\bar{x}|^2 \to 0$, then 
\begin{equation}  \label{eq:happr}
h\big(c(x_r) + \nabla c(x_r) d(x_r,\mu_r)\big) \to h(\bar c).
\end{equation}

\item[{\rm (c)}] When $h$ is convex and lower semicontinuous, the results of
  parts {\rm (a)} and {\rm (b)} hold with $\bar\mu = 0$.
\end{itemize}
\end{theorem}

The next result is a slightly abbreviated version of \cite[Lemma~6.7]{LewW08a}.

\begin{lemma} \label{lem:convergence}
Suppose Assumption \ref{blanket2} holds.
Then for any sequences $\mu_r>0$ and $x_r \to \bar{x}$ such that
$\mu_r |x_r-\bar{x}| \to 0$, and any corresponding sequence of critical points $d_r$
for the subproblem (\ref{pls}) that satisfy the conditions
\begin{equation} \label{eq:dprops}
d_r = O(|x_r - \bar{x}|) ~~\mbox{and}~~ 
h\big(c(x_r) + \nabla c(x_r) d_r\big) \to h(\bar c),
\end{equation}
there exists a bounded sequence of vectors $v_r$ that satisfy
\begin{subequations} \label{eq:hxmopt.r}
\begin{align} 
\label{eq:hxmopt.r.a}
0   & =  \nabla c(x_r)^* v_r + \mu_r d_r, \\
\label{eq:hxmopt.r.b}
v_r &\in \partial h\big(c(x_r) + \nabla c(x_r) d_r\big).
\end{align}
\end{subequations}
\end{lemma}

\noindent
If we assume in addition that $\mu_r>\bar{\mu}$, where $\bar{\mu}$ is
defined in Theorem~\ref{th:proxstep}, the vectors $d_r :=
d(x_r,\mu_r)$ satisfy the properties (\ref{eq:dprops}) and hence the
results of Lemma~\ref{lem:convergence} apply.

We now prove the main result of this section.
\begin{theorem} \label{th:proxid} Suppose Assumption \ref{blanket2}
  holds, and consider a closed set $G \subset \gph(\partial
  h)$. Consider any sequences of scalars $\mu_r >0$ and points $x_r
  \to \bar{x}$ satisfying the condition $\mu_r |x_r-\bar{x}| \to 0$,
  and let $d_r$ be any corresponding sequence of critical points of
  the subproblem \eqnok{pls} satisfying (\ref{eq:dprops}).  Consider
  any corresponding sequence of vectors $v_r$ satisfying the
  conditions \eqnok{eq:hxmopt.r}, and also
\begin{equation} \label{eq:reveal}
\dist \big( (c(x_r) + \nabla c(x_r) d_r, v_r ), G \big) \to 0.
\end{equation}
Then $G$ is an actively sufficient set at $\bar x$. 
\end{theorem}
\begin{proof}
We apply Corollary \ref{sequence}, with ${\cal G} = \{G\}$ and
$c_r := c(x_r) + \nabla c(x_r) d_r$.
Because of the various properties of $x_r$, $d_r$, $v_r$, and $\mu_r$,
from Theorem~\ref{th:proxstep} and Lemma~\ref{lem:convergence}, we
have the following estimates:
\begin{align*}
|x_r-\bar{x}| & \to 0, \\
v_r & \in \partial h(c_r), \\
|c_r - c(x_r) | = \left| \nabla c(x_r) d_r \right| = O(|d_r|) =
O(|x_r-\bar{x}|) & \to 0, \\
|h(c_r) - h(c(\bar{x}))| & \to 0, \\
\left| \nabla c(x_r)^* v_r \right| = \mu_r |d_r| = \mu_r O(|x_r-\bar{x}|) & \to 0 \\
\dist \big( (c_r,v_r), G \big) & \to 0.
\end{align*}
The result follows.
\end{proof}


Note again that Theorem~\ref{th:proxstep} and
Lemma~\ref{lem:convergence} show that vectors $d_r$ satisfying the
conditions of Theorem~\ref{th:proxid} can be obtained when $\mu_r >
\bar{\mu}$, and that we can take $\bar{\mu}=0$ when $h$ is convex and
lower semicontinuous.

As we have seen, in particular in the case of classical nonlinear programming, we typically have in mind some ``natural'' decomposition of the subdifferential graph $\gph(\partial h)$ into the union of a finite family ${\cal G}$ of closed subsets.  We then somehow generate sequences, $\mu_r$, $x_r$, $d_r$, and $v_r$ of the type specified in the theorem, and thereby try to identify actively sufficient sets in ${\cal G}$, preferring smaller sets since the corresponding restricted criticality system is then more refined.  Since ${\cal G}$ is a finite family, Theorem~\ref{th:proxid} guarantees that we must identify at least one actively sufficient set in this way.  However, we may not identify {\em all} actively sufficient sets $G \in {\cal G}$ in this way.   In other words, a sequence of iterates
generated by the algorithm based on (\ref{pls}) and corresponding
multiplier vectors may ``reveal'' some of the actively sufficient sets but not
others. We illustrate this point with an example based on a degenerate
nonlinear optimization problem in two variables.

\begin{Exa} \label{ex:2circle}
{\rm
Consider the map $c:\R^2 \to \R^3$ defined by
\[
c(x) = \left[ \begin{matrix} -x_1 \\ x_1^2+x_2^2-1 \\ (x_1+1)^2 + x_2^2 - 4 
\end{matrix} \right],
\]
and the function $h:\R^3 \to \bar\R$ defined by
\[
h(c) = 
\begin{cases}
c_1 & \mbox{\rm if $c_2,c_3 \le 0$} \\
+\infty & \mbox{\rm otherwise}.
\end{cases}
\]
Minimizing the composite function $h \circ c$ thus amounts to maximizing $x_1$ over the set in $\R^2$ defined by the constraints $|x| \le 1$ and $|x - (-1,0)^T| \le 2$.  The unique minimizer of $h \circ c$ is the point $\bar{x}=(1,0)^T$, at which $c(\bar{x})
= (-1,0,0)^T$.  The set of multiplier vectors is
\[
\partial h(c(\bar{x}))  \cap N(\nabla (c(\bar{x}))^*) ~=~
\Big\{ \alpha\big(1,\frac{1}{2},0\big)^T + (1-\alpha)\big(1,0,\frac{1}{4}\big)^T :
\alpha \in [0,1] \Big\}.
\] 
One decomposition of $\gph(\partial h)$ is as the union of the 
following four closed sets:
\begin{align*}
G^1 &= \{ (c_1,c_2,c_3,1,0,0) : c_2 \le 0, \, c_3 \le 0 \} \\
G^2 &= \{ (c_1,0,c_3,1,v_2,0) : v_2 \ge 0, \, c_3 \le 0 \} \\
G^3 &= \{ (c_1,c_2,0,1,0,v_3) : c_2 \le 0, \, v_3 \ge 0 \} \\
G^4 &= \{ (c_1,0,0,1,v_2,v_3) : v_2 \ge 0, \, v_3 \ge 0 \}.
\end{align*}
(We can refine further, but this suffices for our present purpose.)
In this decomposition, the actively sufficient subsets are $G^2,G^3,G^4$.

The subproblem (\ref{pls}), applied from some point $x=(x_1,0)^T$ with
$x_1$ close to $1$, reduces to
\begin{align*}
\mbox{\rm minimize} \, & -d_1 + \frac{\mu}{2} (d_1^2+d_2^2)  \\
\mbox{\rm subject to} \;\; d_1 & \le \frac{1}{2x_1} - \frac{x_1}{2}, \\
d_1 & \le \frac{2}{x_1+1} - \frac{x_1+1}{2}, \\
d &\in \R^2.
\end{align*}
If $x_1 = 1-\epsilon$ for some small $\epsilon$ (not
necessarily positive), the constraints reduce to
\begin{align*}
d_1 &\le \epsilon + \frac{1}{2} \epsilon^2 + O(\epsilon^3), \\
d_1 &\le \epsilon + \frac{1}{4} \epsilon^2 + O(\epsilon^3).
\end{align*}
Providing $\epsilon \ll \frac{1}{\mu}$, the solution of
the subproblem has $d_1 \approx \epsilon + \frac{1}{4} \epsilon^2$
and $d_2=0$. The corresponding linearized values of $c_2$ and $c_3$ are
\[
c_2(x) + \nabla c_2(x)^Td \approx \frac{1}{4} \epsilon^2, \qquad
c_3(x) + \nabla c_3(x)^Td = 0,
\]
and the corresponding multiplier vector is $v \approx (1,0,\frac{1}{4})^T$. Thus this iterate
``reveals'' the actively sufficient sets $G^3$ and $G^4$, but not $G^2$.

Subsequent iterates generated by this scheme have the identical form
$(1-\epsilon,0)^T$ with successively smaller values of $\epsilon$, so
the sequence satisfies the property (\ref{eq:reveal}) only for $G = G^3$ and $G = G^4$, but not for $G = G^2$.
}
\end{Exa}

Consider again the nonlinear programming formulation of
Section~\ref{sec:classical}.  In that framework, for a given point $x \in \R^n$,  the proximal subproblem \eqnok{pls} is the following quadratic program:
\begin{subequations}
\label{nlp.prox}
\begin{align}
\mbox{\rm minimize} \, f(x) + \nabla   f(x)^T d & +  \frac{\mu}{2} |d|^2 \\
\mbox{\rm subject to} \;\; 
p_i(x) + \nabla p_i(x)^Td & = 0 \;\; (i=1,2,\dotsc,s) \\
q_j(x) + \nabla q_j(x)^Td & \le 0 \;\; (j=1,2,\dotsc,t) \\
                        d & \in \R^n.
\end{align}
\end{subequations}
We derive the following corollary as a simple application of
Theorem~\ref{th:proxid}.
\begin{corollary} \label{co:nlp.prox} Consider the nonlinear program
  {\rm (NLP)}, where the functions $f$, $p_i$ ($i=1,2,\dotsc,s$) and $q_j$
  ($j=1,2,\dotsc,t$) are all $\cC^2$ around the critical point
  $\bar{x}$, and suppose that the Mangasarian-Fromovitz constraint qualification, Assumption~\ref{mf}, holds. Consider sequences of scalars
  $\mu_r>0$ and points $x_r \to \bar{x}$ satisfying $\mu_r
  |x_r-\bar{x}| \to 0$, let $d_r$ be the corresponding (unique) solution of
  (\ref{nlp.prox}), and consider an additional sequence of positive
  tolerances $\epsilon_r \to 0$.  Then for all sufficiently large $r>\bar{r}$, the index set $J(r) \subset \{1,2,\dotsc,t\}$ defined by
\begin{equation} \label{def.Jr}
J(r) ~:=~ 
\big\{ j : 
q_j(x_r) + \nabla q_j(x_r)^T d_r \ge -\epsilon_r  \big\}
\end{equation}
is sufficient for $\bar x$.
\end{corollary}
\begin{proof}
  Suppose the result fails, so that by taking a subsequence, we can
  assume that $J(r)$ is constant: $J(r) = J$ for all $r$, where $J$ is
  not sufficient for $\bar{x}$.  Noting convexity of the function $h$
  defined in Section~\ref{sec:classical} and the equivalence of the
  transversality condition \eqnok{cq} and Assumption~\ref{mf}, we have
  from Theorem~\ref{th:proxid} that the unique solution of the
  subproblem \eqnok{nlp.prox} satisfies $d_r = O(|x_r-\bar{x}|)$ and
\[
h(c(x_r)+\nabla c(x_r) d_r) 
~=~ f(x_r) + \nabla f(x_r)^T d_r 
~\to~ f(\bar{x}) ~=~ h\big(c(\bar{x})\big).
\]
The distance between the point 
\[
\Big( (f(x_r) + \nabla
f(x_r)^T d_r, p(x_r) + \nabla p(x_r) d_r, q(x_r) + \nabla q(x_r)d_r),
(1,\lambda_r,\mu_r) \Big),
\]
and the set $G^J$ defined in (\ref{def.GJ}), approaches zero, where
$\lambda_r$ and $\mu_r$ are the multipliers for the linear constraints
in the subproblem (\ref{nlp.prox}).  We conclude from
Theorem~\ref{th:proxid} that $G^J$ is an actively sufficient set at
$\bar{x}$, so that the index set $J$ is sufficient.  This is a
contradiction.
\end{proof}

Similar results hold for a nonsmooth penalty formulation of the nonlinear program
(NLP).  For example, the $\ell_1$-penalty formulation corresponds to the function $h$
defined as follows:
\[
h(u,v,w) ~=~ u + \nu \Big(\sum_{i=1}^s |v_i|  + \sum_{j=1}^t \max (w_j,0) \Big).
\]
The corresponding proximal subproblem \eqnok{pls} at some given point $x \in \R^n$ is as follows:
\begin{eqnarray*}
\lefteqn {\min_{d \in \R^n}  
f(x) + \nabla f(x)^Td + \frac{\mu}{2} |d|^2 +} \\
& & \mbox{} \hspace{2cm}
\nu \Big( \sum_{i=1}^s \big|p_i(x) + \nabla p_i(x)^Td\big| 
 + \sum_{j=1}^s \max \big(q_j(x)+\nabla q_j(x)^Td,0 \big) \Big),
\end{eqnarray*}
for a given penalty parameter $\nu > 0$.
A result similar to Corollary~\ref{co:nlp.prox} for this formulation
would lead to an identification result like Theorem~3.2 of
\cite{ObeW05a}, provided that $\nu$ is large enough to bound the
$\ell_{\infty}$ norm of all multipliers that satisfy the stationarity
conditions for (NLP). A notable difference, however, is that
\cite[Theorem~3.2]{ObeW05a} uses a trust region of the form 
$\| d \|_{\infty} \le \Delta$ to restrict the size of the solution $d$, whereas
this subproblem uses the prox term $\frac{\mu}{2} |d|^2$. Although the use of
an $\ell_{\infty}$ trust-region allows the subproblem to be formulated
as a linear program, the radius $\Delta$ must satisfy certain
conditions, not easily verified, for the identification result to
hold. By contrast, there are no requirements on $\mu$ in
the subproblems above, beyond positivity.

A possible extension we do not pursue here allows an extra term 
$\frac{1}{2} \ip{d}{Bd}$ for some monotone operator $B$, in addition to the prox
term $\frac{\mu}{2} |d|^2$. This generalization allows SQP type subproblems to be
considered, potentially useful in analyzing algorithms combining
identification and second-order steps into a single iteration
(as happens with traditional SQP methods). 


\end{document}